\documentclass[12pt]{elsarticle}

\usepackage{graphicx,color,enumerate,amssymb}
\usepackage{cauchy}
\usepackage{amsmath}
\usepackage{amssymb}
\usepackage{epsfig}

\makeatletter
\@addtoreset{equation}{section}
\renewcommand\thetable{\@arabic\c@table}
\renewcommand\thefigure{\@arabic\c@figure}
\long\def\@makecaption#1#2{%
  \vskip\abovecaptionskip
  \begin{center}%
  \sbox\@tempboxa{#1: #2}%
  \ifdim \wd\@tempboxa >\hsize
    #1: #2\par
  \else
    \global \@minipagefalse
    \hb@xt@\hsize{\hfil\box\@tempboxa\hfil}%
  \fi
  \end{center}%
  \vskip\belowcaptionskip}
\def\N{{\rm I\kern-.15em N}}
\def\R{{\rm I\kern-.2em R}}
\def\Z{{\rm Z\kern-.26em Z}}
\makeatother
\theoremstyle{plain}
\newtheorem{thm}{Theorem}[section]

\theorembodyfont{\rmfamily}
\newtheorem{rem}[thm]{Remark}
\newtheorem{cor}[thm]{Corollary}

\newcommand{\be}{\begin{eqnarray}}
\newcommand{\ee}{\end{eqnarray}}
\newcommand{\bq}{\begin{eqnarray*}}
\newcommand{\eq}{\end{eqnarray*}}

\def \argmax{\mathop{\hbox{\rm arg max}}}

\newcommand{\bewend}{\hspace*{2mm}\rule{3mm}{3mm}}

\renewcommand{\BBAA}{and}

\newcommand{\RR}{\mathbb{R}}
\newcommand{\PP}{\mathbb{P}}

\newcommand{\BE}{\mathbb{E}}
\newcommand{\BV}{\mathbb{V}}

\newcommand{\LL}{{\textrm L}^2}
\newcommand{\vertk}{\stackrel{{\cal D}}{\longrightarrow}}

\usepackage{color}

\begin{document}
%%%%%%%%%%%%%%%%%%%%%%%%%%%%%%%%%%%%%%%%%%%%%%%%%%%%%%%%%%%%%%%%%%%%%%%%%%%
\begin{center}
{ \LARGE\sc A new class of tests for multinormality with i.i.d. and garch data based on the empirical moment generating function\\
} \vspace*{0.5cm} {\large\sc Norbert Henze}$^{1}$, {\large\sc Mar\'ia Dolores Jim\'enez--Gamero}$^{2}$,\\ \vspace*{0.5cm}
{\it $^{1}$Institute of Stochastics, Karlsruhe Institute of Technology, Karlsruhe, Germany} \\
{\it $^{2}$Department of Statistics and Operations Research, University of Seville, Seville, Spain}\\
%\today
\end{center}

{\small {\bf Abstract.} We generalize a recent class of tests for univariate normality that are based on the empirical moment generating function
to the multivariate setting, thus obtaining a class of affine invariant, consistent and easy-to-use goodness-of-fit tests for multinormality.
The test statistics are suitably weighted $L^2$-statistics, and we provide their asymptotic behavior both for i.i.d.
 observations as well as in the context of testing that the innovation distribution of a multivariate GARCH model is Gaussian.
 We study the finite-sample behavior of the new tests, compare the criteria with alternative existing procedures, and apply the
new procedure to a data set of monthly log returns.}\\ %, and we report applications on real data.}\\ \vspace*{0.3cm}
\vspace*{0.3cm} {AMS
2000 classification numbers:} 62H15, 62G20.
\\ \\
%\newpage

%%%%%%%%%%%%%%%%%%%%%%%%%%%%%%%%%%%%%%%%%%%%%%%%%%%%%%%%%%%%%%%%%%%%%%%%%%%
\section{Introduction}\label{sec_intro}
As evidenced by the papers Arcones (2007), Batsidis et al. (2013), Cardoso de Oliveira and Ferreira (2010), Ebner (2012), Enomoto et al. (2012), Farrel et al. (2007),
Hanusz and Tarasi\'{n}ska (2008, 2012), Henze et al. (2017), Joenssen and Vogel (2014), J\"onsson (2011), Kim (2016), Koizumi et al. (2014), Mecklin and Mundfrom (2005), Pudelko (2005),
Sz\'{e}kely and Rizzo (2005), Tenreiro (2011, 2017), Thulin (2014), Vi\-lla\-se\~nor-Alva and Estrada (2009), Voinov et al. (2016), Yanada et al. (2015),
and Zhou and Shao (2014), there is an ongoing interest in the problem of testing for multivariate normality. Without claiming to be exhaustive, the above list probably covers most
of the publications in this field since the review paper Henze (2002).

Recently, Henze and Koch (2017) provided the lacking theory for a test for {\em univariate} normality suggested by Zghoul (2010).
The purpose of this paper is twofold. First, we generalize the results of Henze and Koch (2017) to the multivariate case, thus obtaining
a class of affine invariant and consistent tests for multivariate normality. Se\-cond, in contrast to that paper (and most of the other publications),
which considered only independent and identically distributed (i.i.d.) observations, we also provide the
asymptotics of our test statistics  in the context of GARCH-type dependence.

To be more specific,
let (for the time being) $X, X_1,X_2, \ldots $ be a sequence of i.i.d. $d$-variate random column vectors that are defined on a common probability space
$(\Omega,{\cal A},\PP)$. We assume that the distribution $\PP^{X}$ is absolutely continuous with respect to Lebesgue measure.
Let N$_d(\mu,\Sigma)$ denote the $d$-variate normal distribution with mean vector $\mu$ and non-degenerate covariance matrix $\Sigma$,
and write ${\cal N}_d$ for the class of all non-degenerate $d$-dimensional normal distributions. A test for multivariate
normality is a test of the null hypothesis
\[
H_0: \ \PP^X  \in {\cal{N}}_d,
\]
and usually such a test should be consistent against any fixed non-normal alternative distribution. Since the class ${\cal N}_d$
is closed with respect to full rank affine transformations, any genuine test statistic $T_n = T_n(X_1,\ldots,X_n)$ based on $X_1,\ldots,X_n$
should also be affine invariant, i.e., we should have
$
T_n(AX_1+b, \ldots, AX_n+b)  =  T_n(X_1,\ldots,X_n)
$
for each nonsingular $d \times d$-matrix $A$ and each $b \in \RR^d$, see Henze (2002) for a critical account on affine invariant tests
for multivariate normality.

In what follows, let $\overline{X}_n = n^{-1}\sum_{j=1}^n X_j$, $S_n = n^{-1}\sum_{j=1}^n (X_j-\overline{X}_n)(X_j-\overline{X}_n)^\top$
denote the sample mean and the sample covariance matrix of $X_1,\ldots,X_n$, respectively, where $\top$  means  transposition of vectors
and matrices. Furthermore, let
\[
 Y_{n,j} = S_n^{-1/2}(X_j - \overline{X}_n), \qquad j=1,\ldots, n,
\]
be the so-called {\em scaled residuals} of $X_1,\ldots,X_n$, which provide an empirical standardization of $X_1,\ldots,X_n$.
 Here, $S_n^{-1/2}$ denotes the unique symmetric square root of $S_n$. Notice that $S_n$ is invertible with probability one provided that
 $n \ge d+1$, see Eaton and Perlman (1973). The latter condition is tacitly assumed to hold in what follows. Letting
 \begin{equation}\label{EMGF}
 M_n(t) = \frac{1}{n} \sum_{j=1}^n \exp\left(t^\top Y_{n,j}\right), \quad t\in \RR^d,
 \end{equation}
 denote the empirical moment generating function of $Y_{n,1}, \ldots,Y_{n,n}$, $M_n(t)$ should be close to
 \[
 m(t)=\exp(\|t\|^2/2),
 \]
  which is the
 moment generating function of the standard normal distribution N$_d(0,\textrm{I}_d)$. Here and in the sequel, $\| \cdot \|$ stands for the Euclidean norm
 on $\RR^d$, and I$_d$ is the unit matrix of order $d$.

The statistic proposed in this paper is the weighted $L^2$-statistic
\begin{equation}\label{teststat}
%T_{n,\beta} = n \int_{\RR^d} \left(M_n(t) - \exp\left( \frac{\|t\|^2}{2}\right) \right)^2 \, w_\beta(t) \, \textrm{d}t,
T_{n,\beta} = n \int_{\RR^d} \left(M_n(t) - m(t) \right)^2 \, w_\beta(t) \, \textrm{d}t,
\end{equation}
where
\begin{equation}\label{weightbeta}
w_\beta(t) = \exp\left(-\beta \|t\|^2 \right),
\end{equation}
and $\beta >1$ is some fixed parameter, the role of which will be discussed later. Notice that $T_{n,\beta}$ is the
'moment generating function analogue' to the BHEP-statistics for testing for multivariate normality (see, e.g., Baringhaus and Henze (1988), Henze and Zirkler (1990), Henze and Wagner (1997)).
The latter statistics originate if one replaces $M_n(t)$ with the empirical {\em characteristic} function of the scaled residuals and $m(t)$ %$\exp(\|t\|^2/2)$
with the characteristic function $\exp(-\|t\|^2/2)$ of the standard normal distribution N$_d(0,{\textrm I}_d)$. For a general account on weighted $L^2$-statistics see, e.g., Baringhaus et al. (2017).

In principle, one could replace $w_\beta$ in (\ref{weightbeta}) with a more general weight function satisfying some general conditions.
The above special choice, however,   leads to a test criterion with certain  extremely appealing features, since straightforward calculations
yield the representation
\begin{eqnarray}\label{Tn}
T_{n,\beta} & = &  \pi^{d/2} \left(\frac{1}{n} \sum_{i,j=1}^n \frac{1}{\beta^{d/2}} \exp\left(\frac{\|Y_{n,i}+Y_{n,j}\|^2}{4 \beta} \right) + \frac{n}{(\beta-1)^{d/2}} \right. \\ \nonumber
&  & \qquad \qquad \left. -2 \sum_{j=1}^n \frac{1}{(\beta -1/2)^{d/2}} \exp\left(\frac{\|Y_{n,j}\|^2}{4 \beta -2} \right)\right),
\end{eqnarray}
which is amenable to computational purposes. Notice that the condition $\beta>1$  is necessary for the integral in (\ref{teststat}) to be finite. Later,
we have to impose the further restriction $\beta >2$ to prove that $T_{n\beta}$ has a non-degenerate limit null distribution as $n \to \infty$.
We remark that $T_{n,\beta}$ is affine invariant since it only depends on the Mahanalobis angles and distances
$Y_{n,i}^\top Y_{n,j}$, $1 \le i,j \le n$.  Rejection of $H_0$ is for large values of $T_{n,\beta}$.

The rest of the paper unfolds as follows. The next section shows that letting $\beta$ tend to infinity in (\ref{teststat})
yields a linear combination of two well-known measures of multivariate skewness. In Section \ref{sec_iid} we derive the limit null distribution
 of $T_{n,\beta}$ in the i.i.d. setting. Section \ref{sec_consist} addresses the question of consistency of the new
 tests against general alternatives, while Section \ref{sec_garch} considers the new criterion in the context of multivariate GARCH models in order to test for normality of innovations,
        and it provides the pertaining large-sample theory. Section \ref{sec_monte} presents a Monte Carlo study that compares the new tests with competing ones, and it considers a real data set
        from the financial market.
         The article concludes with discussions in Section \ref{secconclusion}.

%
%
%
%
%%%%%%%%%%%%%%%%%%%%%%%%%%%%%%%%%%%%%%%%%%%%%%%%%%%%%%%%%%%%%%%%%%%%%%%%%%%%%%%%%%%%%%%%%%%%%%%%%%%%%%%%%%%%%%%%%%%%%%%%%%%%%%%%%%%%%%%%%%%%%%%%%%%%%%%%%%%%%%%%%%%%%%%%%%%%%%%%%%%%%%%%%%%%%%%%%%%%%%%%%%%%%%%%%%%%%%%%

\section{The case $\beta \to \infty$}\label{secbetinfty}
In this section, we show that the statistic $T_{n,\beta}$, after a suitable scaling, approaches a linear combination of two well-known
measures of multivariate skewness  as $\beta \to \infty$.

\begin{thm}\label{beta}
We have
\[
\lim_{\beta \to \infty} \beta^{3+d/2} \, \frac{96 T_{n,\beta}}{n \pi^{d/2}} = 2b_{1,d} +  3\widetilde{b}_{1,d},
\]
where
\[
b_{1,d} = \frac{1}{n^2} \sum_{j,k=1}^n \left(Y_{n,j}^\top Y_{n,k}\right)^3, \qquad \widetilde{b}_{1,d} =  \frac{1}{n^2} \sum_{j,k=1}^n Y_{n,j}^\top Y_{n,k} \, \|Y_{n,j}\|^2 \, \|Y_{n,k}\|^2
\]
are multivariate sample skewness in the sense of Mardia (1970) and
M\'{o}ri, Rohatgi and Sz\'{e}kely (1993), respectively.
\end{thm}
\noindent {\sc Proof}. Let $b_{2,d} = n^{-1}\sum_{j=1}^n \|Y_{n,j}\|^4$ denote multivariate sample kurtosis in the sense of Mardia (1970).
 From (\ref{Tn}) and
\[
\exp(y) = 1 + y + \frac{y^2}{2} + \frac{y^3}{6} + O(y^4)
\]
as $y \to 0$, the result follows by very tedious but straightforward calculations, using the relations
$\sum_{j=1}^n Y_{n,j}  =  0$, $\sum_{j=1}^n \|Y_{n,j}\|^2 =  nd$, $\sum_{j,k=1}^n \|Y_{n,j}+ Y_{n,k}\|^2  =  2n^2d$,
\begin{eqnarray*}
\sum_{j,k=1}^n \|Y_{n,j}+ Y_{n,k}\|^4 \! & \! = \! & \! 2n^2 \left(b_{2,d} + d^2 + 2d\right), \\
\sum_{j,k=1}^n \|Y_{n,j}+ Y_{n,k}\|^4 Y_{n,j}^\top Y_{n,k}  \! & \! = \! & \!  8 n^2 b_{2,d} + 4n^2 b_{1,d} + 2 n^2 \widetilde{b}_{1,d},\\
\sum_{j,k=1}^n \|Y_{n,j}+ Y_{n,k}\|^6 \! & \! = \! & \! 2n \sum_{j=1}^n\!  \|Y_{n,j}\|^6 + 6(d\! +\! 4)n^2 b_{2,d} +8n^2 b_{1,d} + 12 n^2 \widetilde{b}_{1,d}.
\end{eqnarray*}
For the derivation of the second but last expression, see the proof of Theorem 4.1 of Henze et al. (2017).
We stress that although $b_{2,d}$ and $\sum_{j=1}^n \|Y_{n,j}\|^6$ show up in some of the equations above, these terms cancel out in the derivation of the final result. \bewend

\begin{rem} \label{remark.betainf}
Interestingly, $T_{n,\beta}$ exhibits the same limit behavior as $\beta \to\infty$  as both the statistic studied by Henze et al. (2017), which is based
on a weighted $L^2$-distance involving both the empirical characteristic function and the empirical moment ge\-ne\-ra\-ting function, and the BHEP-statistic for testing for multivariate
normality, which is based on the empirical characteristic function, see Theorem 2.1 of Henze (1997). At first sight, Theorem \ref{beta} seems to differ from
Theorem 4 of Henze and Koch (2017) which covers the special case $d=1$,
but a careful analysis shows that -- with the notation $\tau(\beta)$ in that paper -- we have $\lim_{\beta \to \infty} \beta^{7/2} \tau(\beta) =0$.
\end{rem}

%%%%%%%%%%%%%%%%%%%%%%%%%%%%%%%%%%%%%%%%%%%%%%%%%%%%%%%%%%%%%%%%%%%%%%%%%%%%%%%%%%%%%%%%%%%%%%%%%%%%%%%%%%%%%%%%%%%%%%%%%%%%%%%%%%%%%%%%%%%%%%%%%%%%%%%%%%%%%%%%%%%%%%%%%%%%%%%%%%%%%%%%%%%
%
%
%
%
%
%%%%%%%%%%%%%%%%%%%%%%%%%%%%%%%%%%%%%%%%%%%%%%%%%%%%%%%%%%%%%%%%%%%%%%%%%%%%%%%%%%%%%%%%%%%%%%%%%%%%%%%%%%%%%%%%%%%%%%%%%%%%%%%%%%%%%%%%%%%%%%%%%%%%%%%%%%%%%%%%%%%%%%%%%%%%%%%%%%%%%%%%%%%%%%%

\section{Asymptotic null distribution in the i.i.d. case}\label{sec_iid}
In this section we consider the case that $X_1,X_2, \ldots $ are i.i.d. $d$-dimensional random
vectors with some non-degenerate normal distribution. The key observation for deriving the limit distribution of $T_{n,\beta}$
is the fact that
\[%\label{unl2}
T_{n,\beta} \ = \ \int_{\RR^d} W_n(t)^2 \, w_\beta(t) \, \textrm{d}t,
\]
where
\begin{equation}\label{defwn}
%W_n(t) = \sqrt{n}\left(M_n(t) - \exp\left(\|t\|^2/2\right)\right), \quad t \in \RR^d,
W_n(t) = \sqrt{n}\left(M_n(t) - m(t)\right), \quad t \in \RR^d,
\end{equation}
with $M_n(t)$ given in (\ref{EMGF}).
Notice that $W_n$ is a random element of the  Hilbert space
\begin{equation}\label{hildef}
\textrm{L}^2_\beta := \textrm{L}^2(\RR^d,{\cal B}^d,w_\beta(t)\textrm{d}t)
\end{equation}
of (equivalence classes of) measurable functions $f:\RR^d \rightarrow \RR$
that are square integrable with respect to the finite measure on the $\sigma$-field ${\cal B}^d$ of Borel sets of $\RR^d$ given by the weight function $w_\beta$
defined in (\ref{weightbeta}).
The resulting norm in $\textrm{L}_\beta^2$ will be denoted by
$\|f\|_{\textrm{L}_\beta^2} = \sqrt{\langle f,f \rangle}$. With this notation, $T_{n,\beta}$ takes the form
\begin{equation}\label{unl2}
T_{n,\beta} \ = \ \|W_n\|^2_{\textrm{L}^2_\beta}.
\end{equation}
Writing "$\vertk$" for convergence in distribution of random vectors and stochastic processes, the main result of this section
is as follows.

\smallskip
\begin{thm}{\rm{(}}Convergence of $W_n$ under $H_0${\rm{)}}\label{mainthm}\\
Suppose that $X$ has some non-degenerate $d$-variate normal distribution, and that $\beta >2$ in (\ref{weightbeta}). Then there is a centred Gaussian random element $W$ of
$\LL$ having covariance kernel
\[
C(s,t)   =   \exp\left(\frac{\|s\|^2+\|t\|^2}{2} \right) \left({\rm{e}}^{s^\top t} -1 - s^\top t - \frac{\left(s^\top t \right)^2}{2} \right), \quad s,t \in \RR^d,
\]
so that $W_n \vertk W$ as $n \to \infty$.
\end{thm}

In view of \eqref{unl2}, the Continuous Mapping Theorem yields the following result.

\begin{cor} \label{cor1}
If $\beta >2$, then, under the null hypothesis $H_0$,
\[
 T_{n, \beta} \vertk \|W\|^2_{\LL_\beta} \ \ \textrm{as} \ n \to \infty.
 \]
%where $W_G$ is the Gaussian random element figuring in Theorem \ref{mainthmG}.
\end{cor}

\begin{rem}
The distribution of $T_{\infty,\beta} := \|W\|^2_{\LL_\beta}$ (say) is that of $\sum_{j=1}^\infty \lambda_jN_j^2$, where $\lambda_1, \lambda_2, \ldots $ are the
positive eigenvalues of the integral operator $f \mapsto Af$  on $\LL_\beta$ associated with the kernel $C$ given in Theorem \ref{mainthm},
i.e., $(Af)(t) \! = \! \int \! C(s,t) f(s) \exp(-\beta \|s\|^2) \textrm{d}s$,
and $N_1,N_2, \ldots $ are i.i.d. standard normal random variables.
We did not succeed in obtaining explicit solutions of this equation.
However, since
\begin{eqnarray*}
\mathbb E (T_{\infty,\beta}) & = & \int_{\RR^d} C(t,t) \, w_\beta(t) \, \textrm{d}t,\\
\BV(T_{\infty,\beta}) & = & 2 \int_{\RR^d} \int_{\RR^d} C^2(s,t) w_\beta(s) w_\beta(t) \, \textrm{d}s \textrm{d}t
\end{eqnarray*}
(see Shorack and Wellner, 1986, p. 213), tedious but straighforward manipulations of integrals yield the following result, which generalizes Theorem 2 of Henze and Koch (2017).
\end{rem}
\begin{thm}\label{thmevar}
If $\beta >2$ we have
\begin{enumerate}
\item[a)]
\[
\BE(T_{\infty,\beta})  = \pi^{d/2} \left(\frac{1}{(\beta -2)^{d/2}} - \frac{1}{(\beta -1)^{d/2}} - \frac{d}{2(\beta -1)^{d/2+1}} - \frac{d(d+2)}{8 (\beta -1)^{d/2+2}}\right),\\
\]
\item[b)]
\begin{eqnarray*}
\BV(T_{\infty,\beta}) & = &   2\pi^d\left(\frac{1}{(\beta(\beta-2))^{d/2}} - \frac{2^{d+1}}{\eta^{d/2}} - \frac{(1+2d)2^d}{\eta^{d/2+1}} - \frac{d(d+2)2^d}{\eta^{d/2 +2}} \right. \\
& & \qquad \left. + \frac{1}{(\beta -1)^d} + \frac{d}{2(\beta -1)^{d+2}} + \frac{3d(d+2)}{64 (\beta -1)^{d+4}}\right),
\end{eqnarray*}
\end{enumerate}
where $\eta= 4(\beta -1)^2-1$.
\end{thm}
%%%%%%%%%%%%%%%%%%%%%%%%%%%%%%%%%%%%%%%%%%%%%%%%%%%%%%%%%%%%%%%%%%%%%%%%%%%%%%%%%%%%%%%%%%%%%%%%%%%%%%%%%%%%%%%%%%%%%%%%%%%%%%%%%%%%%%%%%%%%%%%%%%%%%%%%%%%%%%%%%%%%%%%%%%%%%%%%%%%%%%%%%%%%%%%%%%%%%%%%%%%%%%%%%
%
%
%
%
%
%%%%%%%%%%%%%%%%%%%%%%%%%%%%%%%%%%%%%%%%%%%%%%%%%%%%%%%%%%%%%%%%%%%%%%%%%%%%%%%%%%%%%%%%%%%%%%%%%%%%%%%%%%%%%%%%%%%%%%%%%%%%%%%%%%%%%%%%%%%%%%%%%%%%%%%%%%%%%%%%%%%%%%%%%%%%%%%%%%%%%%%%%%%%%%%%%%%%%%%%%%%%%%%%%%%%%

\vspace*{4mm}

\noindent {\sc Proof} of Theorem \ref{mainthm}. In view of affine invariance, we assume w.l.o.g. that the distribution of $X$ is N$_d(0,{\rm{I}}_d)$. In Henze et al. (2017), the authors considered  the
``exponentially down-weighted empirical moment generating function process''
\begin{equation}\label{defan}
A_n(t) = \exp\left(- \frac{\|t\|^2}{2}\right) \, M_n(t), \quad t \in \RR^d.
\end{equation}
Notice that, with the notation given in (\ref{hildef}), we have
\[
\|A_n\|^2_{\textrm{L}^2_\beta} = \|M_n\|^2_{\textrm{L}^2_{\gamma}},
\]
where $\gamma = \beta -1$
From display (10.5) and Propositions 10.3 and 10.4  of Henze et al. (2017) we have
\[
%A_n(t) = \exp\left(- \frac{\|t\|^2}{2} \right) \sqrt{n} \left(\frac{1}{n}\sum_{j=1}^n {\rm{e}}^{t^\top X_j} - \exp\left(\frac{\|t\|^2}{2} %\right)\right) + V_{n}(t) + R_n(t),
A_n(t) = \exp\left(- \frac{\|t\|^2}{2} \right) \sqrt{n} \left(\frac{1}{n}\sum_{j=1}^n {\rm{e}}^{t^\top X_j} - m(t)\right) + V_{n}(t) + R_n(t),
\]
where $\int_{\RR^d} R_n^2(t) w_\gamma (t)\textrm{d} t = o_\PP(1)$,
and
\[
V_{n}(t) = - \frac{1}{2\sqrt{n}} \sum_{j=1}^n \left((t^\top X_j)^2 - \|t\|^2\right) - \frac{1}{\sqrt{n}} \sum_{j=1}^n t^\top X_j.
\]
%Putting $m(t) = \exp(\|t\|^2/2)$, (\ref{defan}) and the representation of $A_n$ as a sum yield
Display (\ref{defan}) and the representation of $A_n$ as a sum yield
\[
W_n(t) = \frac{1}{\sqrt{n}} \sum_{j=1}^n Z_j(t) + m(t)R_n(t),
\]
where
\[
Z_j(t) = {\rm{e}}^{t^\top X_j} - m(t) - \frac{m(t)}{2} \left( (t^\top X_j)^2 - \|t\|^2\right) - m(t) t^\top X_j.
\]
Notice that $Z_1,Z_2, \ldots $ are i.i.d. centred random elements of
L$^2_\beta$. Since
\[
\int_{\RR^d} (m(t)R_n(t))^2 w_\beta(t)\, \textrm{d}t = \int_{\RR^d} R_n^2(t) w_\gamma(t) \, \textrm{d} t = o_\PP(1),
\]
 a Central Limit Theorem  in Hilbert spaces (see e.g., Bosq, 2000) shows that there is a centered Gaussian random element $W$
 of L$^2_\beta$, so that $W_n \vertk W.$  Using the fact that
 $t^\top X$ has the normal distribution N$(0,\|t\|^2)$ and the relations
 \begin{eqnarray*}
 \BE \left[{\rm{e}}^{s^\top X} (t^\top X)^2 \right] & = & m(s)\left( (s^\top t)^2 + \|t\|^2\right),\\
 \BE \left[ {\rm{e}}^{s^\top X} t^\top X \right] & = & m(s) s^\top t,\\
 \BE\left[ (s^\top X)^2 (t^\top X)^2 \right] & = & 2 (s^\top t)^2 + \|s\|^2 \, \|t\|^2,
 \end{eqnarray*}
some straightforward algebra shows that the covariance kernel $C(s,t)$ figuring in the statement of Theorem
 \ref{mainthm} equals $\BE Z_1(s)Z_1(t)$. \ $ \bewend$

\section{Consistency}\label{sec_consist}

The next result shows that the test for multivariate normality based on $T_{n,\beta}$ is consistent against general alternatives.

\begin{thm}\label{thmconsist}
Suppose $X$ has some absolutely continuous distribution, and that $M_X(t) = \BE [\exp(t^\top X)] < \infty$, $t \in \RR^d$.
Furthermore, let $\widetilde{X} = \Sigma^{-1/2}(X-\mu)$, where $\mu = \BE(X)$  and $\Sigma^{-1/2}$ is the symmetric square
root of the inverse of the covariance matrix $\Sigma$ of $X$. Letting
$M_{\widetilde{X}}(t) = \BE [\exp(t^\top \widetilde{X})]$, % and again putting $m(t) = \exp(\|t\|^2/2)$,
we have
\[ %\begin{equation}\label{consistin}
\liminf_{n \to \infty} \frac{T_{n,\beta}}{n} \ \ge \  \int_{\RR^d} \left(M_{\widetilde{X}}(t) -m(t)  \right)^2 \, w_\beta(t) \, \rm{d} t
\] %\end{equation}
almost surely.
\end{thm}

\noindent {\sc Proof}.  Because of affine invariance we may w.l.o.g. assume $\BE X =0$ and $\Sigma = \textrm{I}_d$.
Fix $K>0$ and put  $M_n^\circ(t)  = n^{-1}\sum_{j=1}^n \exp(t^\top X_j)$. From the proof of Theorem 6.1 of Henze et al. (2017) we have
\[
\lim_{n\to \infty} \max_{\|t\| \le K} \big{|} M_n(t) - M_n^\circ(t) \big{|} = 0
\]
$\PP$-almost surely. Now, the strong law of large numbers in the Banach space of continuous functions on $B(K):=\{t \in \RR^d:\|t\| \le K\}$ and Fatou's lemma yield
\begin{eqnarray*}
\liminf_{n\to \infty} \frac{T_{n,\beta}}{n} & \ge & \liminf_{n \to \infty} \int_{B(K)} \left(M_n(t) - m(t)\right)^2 w_\beta(t) \, \textrm{d}t\\
& \ge & \int_{B(K)} \left( \BE \textrm{e}^{t^\top X}- m(t)\right)^2 w_\beta(t) \, \textrm{d}t
\end{eqnarray*}
$\PP$-almost surely. Since $K$ is arbitrary, the assertion follows. \ $\bewend$

\vspace*{5mm}

Now, suppose that $X$ has an alternative distribution (which is assumed to be standardized)
satisfying the conditions of Theorem \ref{thmconsist}. Since $\BE \exp(t^\top X) - m(t) \neq 0$ for at least one $t$,
Theorem \ref{thmconsist} shows that $\lim_{n\to \infty} T_{n,\beta} = \infty$ $\PP$-almost surely. Since, for any given nominal level $\alpha \in (0,1)$,
the sequence of critical values of a level-$\alpha$-test based on $T_{n\beta}$ that rejects $H_0$ for large values of $T_{n,\beta}$ converges according to
Theorem \ref{mainthm}, this test is consistent against such an alternative. It should be 'all the more consistent' against any distribution not satisfying the
 conditions of Theorem \ref{thmconsist} but, in view of the reasoning given in Cs\"org\H{o} (1986), the behavior of $T_{n,\beta}$ against such alternatives  is a difficult problem.

%%%%%%%%%%%%%%%%%%%%%%%%%%%%%%%%%%%%%%%%%%%%%%%%%%%%%%%%%%%%%%%%%%%%%%%%%%%%%%%%%%%%%%%%%%%%%%%%%%%%%%%%%%%%%%%%%%%%%%%%%%%%%%%%%%%%%%%%%%%%%%%%%%%%%%%%%%%%%%%%%%%%%%%%%%%%%%%%%%%%%%%%%%%%%%%%%%%%%%%%%%%%%%%%%%
%
%
%
%
%
%%%%%%%%%%%%%%%%%%%%%%%%%%%%%%%%%%%%%%%%%%%%%%%%%%%%%%%%%%%%%%%%%%%%%%%%%%%%%%%%%%%%%%%%%%%%%%%%%%%%%%%%%%%%%%%%%%%%%%%%%%%%%%%%%%%%%%%%%%%%%%%%%%%%%%%%%%%%%%%%%%%%%%%%%%%%%%%%%%%%%%%%%%%%%%%%%%%%%%%%%%%%%%%%%%%%

\section{Testing for normality in GARCH models}\label{sec_garch}
In this section we consider the multivariate GARCH (MGARCH) model
\begin{equation} \label{GARCH}
X_j=\Sigma_j^{1/2}(\theta)\varepsilon_j, \quad j\in \mathbb{Z},
\end{equation}
where $\theta \in \Theta \subseteq \mathbb R^v$ is a $v$-dimensional vector of  unknown parameters.
The unobservable random errors or innovations $\{\varepsilon_j, \,j \in \mathbb{Z}\}$  are i.i.d. copies
 of a $d$-dimensional random vector $\varepsilon$, which is assumed to have mean zero and unit covariance matrix. Hence \[\Sigma_j(\theta)=\Sigma(\theta; X_{j-1},X_{j-2}, \ldots)\]
 is the conditional variance of $X_j$, given $X_{j-1},X_{j-2}, \ldots$. The explicit expression of $\Sigma_j(\theta)$  depends on the assumed MGARCH model (see, e.g., Francq and Zako{\"i}an, 2010, for a detailed description of several relevant models).
The interest in testing for normality of the innovations stems from the fact that this distributional assumption is made in some applications, and that,
 if  erroneously accepted, some  inferential procedures can lead to wrong conclusions (see, e.g., Spierdijk, 2016, for the effect on the assessment of standard risk measures such as the value at risk).

 Therefore, an important step
in the analysis of GARCH models is to check whether the data support  the distributional hypotheses made
on the innovations. Because of this reason, a number of goodness-of-fit tests have been proposed
for the innovation distribution.
The papers by Klar et al. (2012) and Ghoudi and R\'emillard (2014) contain an extensive review of such
tests as well as some numerical comparisons between them for the special case of testing for univariate normality. The proposals for testing  goodness-of-fit in the multivariate case are rather scarce.

The class of GARCH models has been proved to be particularly valuable in mo\-de\-ling financial
data. As discussed, among others in Rydberg (2000), one of the stylized features of financial data
is that they are heavy-tailed.  From an extensive simulation study (a summary is reported in Section \ref{sec_monte}), we learnt that, for i.i.d. data, the test of normality
based on $T_{n,\beta}$ exhibits a high power against heavy-tailed distributions. Because of these reasons, this section is devoted to adapt that procedure to testing for normality
 of the innovations based on  data $X_1, \ldots, X_n$ that are  driven by equation (\ref{GARCH}).
Therefore, on the basis of the  observations, we wish to test the null
hypothesis
\[
{H}_{0,G}: \ {\rm{The \ law \ of}} \  \varepsilon  \ {\rm{is}} \ \textrm{N}_d(0,{ {\rm I}}_d).
\]
against general alternatives. Notice that ${{H}}_{0,G}$ is equivalent to the hypothesis that, conditionally on $\{X_{j-1},X_{j-2},\ldots\}$, the law of $X_j$  is N$_d(0,\Sigma_j(\theta))$,
for some $\theta \in \Theta$. Two main differences with respect to the i.i.d. case are: (a) the innovations in (\ref{GARCH}) are assumed to be centered at zero with unit covariance matrix;  and (b) the conditional covariance matrix $\Sigma_j(\theta)$ of $X_j$ is time-varying in a way that depends on the unknown parameter $\theta$ and  on past observations.

Notice that although  ${{H}}_{0,G}$  is about the distribution of $\varepsilon$, the innovations  are unobservable in the context of model (\ref{GARCH}). Hence any inference on the distribution of the innovations should be based on residuals
\begin{equation} \label{residuoss}
{\widetilde {\varepsilon}}_j(\widehat{\theta}_n)={\widetilde{ \Sigma}}_j^{-1/2}(\widehat{\theta}_n)X_j, \quad 1 \leq j \leq n.
\end{equation}
 Recall that $\Sigma_j(\theta)=\Sigma(\theta; X_{j-1},X_{j-2},\dots)$, but we only observe $X_1, \ldots, X_n$. Therefore, to estimate  $\Sigma_j(\theta)$, apart from a suitable estimator $\widehat{\theta}_n$ of $\theta$, we also need to
specify  values for $\{X_j, \ j\leq 0\}$, say $\{\widetilde X_j, \ j\leq 0\}$. So we write
${\widetilde{ \Sigma}}_j({\theta})$ for  $\Sigma(\theta; X_{j-1}, \ldots, X_{1},\widetilde{X}_0,\widetilde{X}_{-1} \ldots )$.
Under certain conditions, these arbitrarily fixed initial values are asymptotically irre\-le\-vant.

Taking into account that the innovations have mean zero and unit covvariance matrix, we will work directly with the residuals, without standardizing them.
Let $M_n^G$ be defined  as $M_n$ in \eqref{EMGF} by replacing $Y_{n,j}$ with ${\widetilde {\varepsilon}}_j(\widehat{\theta}_n)$, $1\leq j \leq n$,  and define $T_{n, \beta}^G$  as $T_{n, \beta}$  in \eqref{unl2} with $W_n$ changed for $W_n^G$, where $W_n^G$ is defined as $W_n$ in \eqref{defwn} with  $M_n$ replaced by $M_n^G$. In order to derive the asymptotic null distribution of $W_n^G$ we will make the assumptions (A.1)--(A.6) below.
  In the sequel, $C>0$ and $\varrho$, $0< \varrho < 1$,  denote generic constants,
 the values of which may vary across the text, $\theta_0$ stands for the true value of $\theta$, and for any  matrix $A=(a_{kj})$, $\| A \|=\sum_{k,j}|a_{k j}|$
 denotes the $l_1$-norm (we use the same notation as for the Euclidean norm of vectors).

 \begin{enumerate}
\item[(A.1)] The estimator $\widehat{\theta}_n$ satisfies
$
\sqrt{n}(\widehat{\theta}_n-\theta_0)= n^{-1/2} \sum_{j=1}^nL_j +o_{\mathbb{P}}(1),%(\theta; \varepsilon_j, \varepsilon_{j-1}, \ldots)+o_{\mathbb{P}}(1),
$
where $L_j=h_j g_j$, $g_j=g(\theta_0; \varepsilon_j)$ is a vector of $d^2$ measurable functions such that $ \mathbb E ( g_j)=0$ and $\mathbb E ( g_j^\top g_j)<\infty$, and
$h_j=h (\theta_0; \varepsilon_{j-1}, \varepsilon_{j-2}\ldots)$ is a $v\times d^2$-matrix of measurable functions satisfying $\mathbb E ( \|h_j h_j^\top \|^2)<\infty$,

\item[(A.2)] $\sup_{\theta\in\Theta}\left\|\widetilde{\Sigma}^{-1/2}_{j}(\theta)\right\|\leq C, \quad \sup_{\theta\in\Theta}\left\|\Sigma^{-1/2}_{j}(\theta)\right\|\leq C\quad \PP\mbox{-a.s.}$,

\item[(A.3)] $\sup_{\theta\in\Theta}\|\Sigma^{1/2}_{j}(\theta)-\widetilde{\Sigma}^{1/2}_{j}(\theta)\| \leq C\varrho^j$,

\item[(A.4)] $\mathbb E \left\|X_j\right\|^\varsigma<\infty$ and $ \mathbb E\left\|\Sigma^{1/2}_{j}(\theta_0)\right\|^\varsigma<\infty$  for some $\varsigma>0$,

\item[(A.5)] for any sequence $x_1,x_2,\dots$ of vectors of $\mathbb{R}^d$, the function $\theta\mapsto	
	\Sigma^{1/2}(\theta; x_1,x_2,\dots)$ admits continuous second-order derivatives,

\item[(A.6)] for some neighborhood $V(\theta_0)$ of  $\theta_0$, there exist $p> 1$, $q> 2$ and $r> 1$ so that $2p^{-1}+2r^{-1}=1$ and $4q^{-1}+2r^{-1}=1$, and
    \begin{eqnarray*}
&& \mathbb E \sup_{\theta\in V(\Theta)}\left\|\sum_{k,\ell=1}^v\Sigma_j^{-1/2}(\theta)\frac{\partial^2\Sigma^{1/2}_j(\theta)}{\partial \theta_k\partial \theta_\ell}\right\|^p<\infty,\\
&& \mathbb E \sup_{\theta\in V(\Theta)}\left\|\sum_{k=1}^v\Sigma_j^{-1/2}(\theta)\frac{\partial\Sigma^{1/2}_j(\theta)}{\partial \theta_k}\right\|^q<\infty,\\
&& \mathbb E \sup_{\theta\in V(\Theta)}\left\|\Sigma_j^{1/2 }(\theta_0)\Sigma_j^{-1/2}(\theta)\right\|^r<\infty.
	\end{eqnarray*}
\end{enumerate}

The next result gives the asymptotic null distribution of $W_n^G$.

\begin{thm}{\rm{(}}Convergence of $W_n^G$ under $H_{0,G}${\rm{)}}\label{mainthmG}\\
Let $\{X_j\}$ be a strictly stationary process satisfying \eqref{GARCH}, with $X_j$ being measurable with respect to the sigma-field generated by $\{\varepsilon_u,u\leq j\}$. Assume  that
    (A.1)--(A.6) hold and that $\beta>2$. Then under the null hypothesis ${{H}}_{0,G}$,  there is a centered Gaussian random element $W_G$ of
$\LL_{\beta}$, having covariance kernel
$
C_G(s,t)   = cov(U(t),U(s)),%
$
 so that $W_n^G \vertk W_G$ as $n \to \infty$, where
\[
%U(t)=\exp(t^\top \varepsilon_1)-\exp(\|t\|^2/2)-\exp(\|t\|^2/2)a(t)^\top L_1,
U(t)=\exp(t^\top \varepsilon_1)-m(t)-m(t)a(t)^\top L_1,
\]
$a(t)^\top =(t^\top \mu_1 t, \ldots, t^\top \mu_v t)$, $\mu_k=\mathbb{E}[A_{1k}(\theta_0)]$, $A_{1k}(\theta)=\Sigma_1^{-1/2}(\theta)\frac{\partial}{\partial \theta_k} \Sigma_1^{1/2}(\theta)$, $1\leq k \leq v$.
\end{thm}

From Theorem \ref{mainthmG} and the Continuous Mapping Theorem we have the following corollary.

\begin{cor} \label{cor2}
Under the assumptions of Theorem \ref{mainthmG}, we have
\[ T_{n, \beta}^G \vertk \|W_G\|^2_{\LL_{\beta}} \ \ \textrm{as} \ n \to \infty.\]
\end{cor}

The standard estimation method for the parameter $\theta$ in  GARCH models is the quasi maximum likelihood estimator (QMLE), defined as
\[ %\be \label{defestimator}
\widehat{\theta}_n=\argmax_{\theta \in \Theta} {L}_n(\theta),
\] %\ee
where
\[
L_n(\theta)=-\frac{1}{2} \sum_{j=1}^n \widetilde{\ell}_j (\theta),\quad \widetilde{\ell}_j(\theta )=X^\top_j
\widetilde{ \Sigma}_j(\theta)^{-1}X_j +\log\left|\widetilde{ \Sigma}_j(\theta)\right|.
\]
Comte and Leiberman (2003) and  Bardet and Wintenberger (2009), among others, have shown that under certain mild regularity conditions the QMLE satisfies (A.1) for general MGARCH models.

As observed before, there are many MGARCH parametrizations for the matrix $\Sigma_j(\theta)$. Nevertheless,  there exist only  partial theoretical results  for such models.
The Constant Conditional Correlation  model, proposed by Bollerslev (1990) and  extended by Jeantheau (1998),   is an exception, since its properties have been thoroughly studied.
This model decomposes the conditional covariance matrix figuring in (\ref{GARCH}) into conditional standard
deviations and a conditional correlation matrix, according to
$
\Sigma_j(\theta_0)={ {D}}_j(\theta_0) R_0 {{D}}_j(\theta_0),
$
where ${{D}}_j(\theta_0)$ and $R_0$ are $d\times d$-matrices,  $R_0$ is a correlation matrix, and  ${{D}}_j(\theta_0)$ is a diagonal
 matrix so that $\sigma^2_j(\theta)=\mbox{diag}\left\{D^2_j(\theta)\right\}$ with
\begin{equation} \label{fla.CCCGARCH}
\sigma^2_j(\theta)={b}+\sum_{k=1}^p {{B}}_k X_{j-k}^{(2)}+\sum_{k=1}^q {{\Gamma}}_k \sigma^2_{j-k}(\theta).
\end{equation}
Here, $X_{j}^{(2)}= X_{j} \odot X_{j}$, where $\odot$ denotes the Hadamard product, that is, the element by element product,
${{b}}$ is a vector of dimension $d$ with positive elements, and $\{B_k\}_{k=1}^p$ and $\{{{\Gamma}}_k\}_{k=1}^q$ are  $d\times d$ matrices
with non-negative elements. This model will be referred to as CCC-GARCH($p$,$q$).
Under certain weak assumptions, the QMLE for the parameters in this model satisfies (A.1), and (A.2)--(A.6) also hold, see Francq and Zako{\"i}an (2010) and Francq et al. (2017).

The asymptotic null distribution of $T_{n, \beta}^G $ depends on   the equation defining the GARCH model and on $\theta_0$ through the quantities $\mu_1, \ldots, \mu_v$,
as well as on which estimator of $\theta$ has been employed. Therefore, the asymptotic null distribution cannot be used to approximate the null distribution of $T_{n, \beta}^G $. Following
Klar et al. (2012), we will estimate the null distribution  of $T_{n, \beta}^G $  by using the following
 parametric bootstrap algorithm:
 \begin{itemize} \itemsep=0pt
\item[(i)] Calculate $\widehat{\theta}_n=\widehat{\theta}_n(X_1, \ldots, X_n)$, the residuals  $\widetilde{\varepsilon}_1,\ldots,\widetilde{\varepsilon}_n$ and the test statistic ${T}_{n,\beta}^G=
{T}_{n,\beta}^G(\widetilde{\varepsilon}_1,\ldots,\widetilde{\varepsilon}_n)$.

\item[(ii)] Generate vectors ${\varepsilon}_1^*,\ldots,{\varepsilon}_n^*$ i.i.d. from a N$_d(0,{ {\rm I}}_d)$ distribution. Let  $X_j^*=\widetilde{\Sigma}_j^{1/2}(\widehat{\theta})\varepsilon_j^*$, $j=1,\ldots, n$.

\item[(iii)]    Calculate $\widehat{\theta}_n^*=\widehat{\theta}_n(X_1^*, \ldots, X_n^*)$, the residuals  $\widetilde{\varepsilon}_1^*,\ldots,\widetilde{\varepsilon}_n^*$, and approximate the null distribution of ${T}^G_{n,\beta}$ by means of the conditional distribution, given the data, of  ${T}^{G*}_{n,\beta}={T}^G_{n,\beta}(\widetilde{\varepsilon}_1^*,\ldots,\widetilde{\varepsilon}_n^*)$.
\end{itemize}

In practice, the approximation in step (iii)   is carried out by generating a large number of bootstrap replications of the test statistic  ${T}_{n,\beta}^{G}$,
whose empirical distribution function is used to estimate the null distribution of ${T}^G_{n,\beta}$.  Similar steps to those given in the proof of Theorem  \ref{mainthmG} show that
  if one assumes that (A.1)--(A.6) continue to hold when  $\theta_0$ is replaced by $\theta_n$, with $\theta_n\to \theta_0$ as $n\to \infty$, and $\varepsilon \sim\textrm{N}_d(0,{ {\rm I}}_d)$, then  the conditional distribution of ${T}^{G*}_{n,\beta}$, given the data, converges in law to $ \|W_G\|^2_{\LL_{\beta}}$, with $W_G$ as defined in Theorem  \ref{mainthmG}. Therefore, the
above bootstrap procedure provides a consistent null distribution estimator.

\begin{rem}
The practical application of the above bootstrap null distribution estimator entails that  the parameter estimator of $\theta$ and the residuals must be calculated for each bootstrap resample, which  results in a time-consuming procedure. Following the approaches in Ghoudi and R\'emillard (2014) and Jim\'enez-Gamero and Pardo-Fern\'andez (2017) for other goodness-of-fit tests for univariate GARCH models, we could use a  weighted bootstrap null distribution estimator in the sense of Burke (2000). From a computational point of view, it provides a more efficient estimator. Nevertheless, it can be verified that the consistency of the weighted bootstrap null distribution estimator of ${T}_{n,\beta}^G$ requires the existence of the moment generating function of the true distribution generating the innovations, which is a rather strong condition, specially taking into account that the alternatives of interest are heavy-tailed.
\end{rem}

As in the i.i.d. case, the next result shows that the test for multivariate normality based on $T_{n,\beta}^G$ is consistent against general alternatives.

\begin{thm}\label{thmconsistG}
Let $\{X_j\}$ be a strictly stationary process satisfying \eqref{GARCH}, with $X_j$ being measurable with respect to the sigma-field generated by $\{\varepsilon_u,u\leq j\}$. Assume  that
    (A.1)--(A.6) hold, that
$\varepsilon$ has some absolutely continuous distribution, and that $M_\varepsilon(t) = \BE [\exp(t^\top \varepsilon)] < \infty$, $t \in \RR^d$. We then have
\[ %\begin{equation}\label{consistin}
\liminf_{n \to \infty} \frac{T_{n,\beta}^G}{n} \ \ge \  \int_{\RR^d} \left(M_{\varepsilon}(t) -m(t)  \right)^2 \, w_\beta(t) \, \rm{d} t
\] %\end{equation}
in probability.
\end{thm}

Similar comments to those made after Theorem \ref{thmconsist} for the i.i.d. case can be done in this setting.

\noindent{\sc Proof} of Theorem \ref{mainthmG}. From the  proof of Theorem 7.1 in Henze et al. (2017), it follows that
$W_n^G(t)=W_{1,n}^G(t)+r_{n,1}(t)$,  with
$W_{1,n}^G(t)=n^{-1/2}\sum_{j=1}^nV_j(t)$,
\[V_j(t)=\exp(t^\top\varepsilon_j)-m(t)a(t)^\top\sqrt{n}(\widehat{\theta}_{n}-\theta_{0})-m(t),\]
and $\|r_{n,1}\|_{\LL_{\beta}}=o_{\mathbb{P}}(1)$.
By Assumption A.1,
$W_{1,n}^G(t)=W_{2,n}^G(t)+r_{n,2}(t)$,  with
$W_{2,n}^G(t)=n^{-1/2}\sum_{j=1}^nU_j(t)$,
\[U_j(t)=\exp(t^\top\varepsilon_j)-\exp(-\|t\|^2/2)a(t)^\top L_j-\exp(-\|t\|^2/2),\]
and $\|r_{n,2}\|_{\LL_{\beta}}=o_{\mathbb{P}}(1)$.

To prove the result we will apply Theorem 4.2 in Billingsley (1968) to  $\{W_{2,n}^G(t), t \in \mathbb{R}^d \}$ by showing that (a) for each positive $M$,  $\{W_{2,n}^G(t), t \in B(K) \}$  converges in law to  $\{W^G(t), t \in B(K) \}$ in $C(B(K) )$, the Banach space of real-valued continuous functions on $B(K):=\{t \in \RR^d:\|t\| \le K\}$, endowed with the supremum norm; (b) for each $\varepsilon>0$, there is a positive $K$ so that
\begin{equation} \label{resto1}
\int_{\mathbb{R}^d\setminus B(K)} \mathbb E \left[W_{2,n}^G(t)^2\right] w_\beta(t) \, \textrm{d}t < \varepsilon,
\end{equation}
\begin{equation} \label{resto2}
\int_{\mathbb{R}^d\setminus B(K)} \mathbb E\left[W^G(t)^2\right] w_\beta(t) \, \textrm{d}t < \varepsilon.
\end{equation}

\noindent \underline{Proof of (a)}: By applying the central limit theorem for martingale differences, the  finite-dimensional distributions of
$\{W_{2,n}^{G}(t), t \in \mathbb{R}^d \}$ converge to those of $\{W_{G}(t), t \in \mathbb{R}^d \}$. Hence, to prove (a) we must show that $\{W_{2,n}^{G}(t), t \in B(K) \}$ is tight.
 With this aim we write $W_{2,n}^{G}(t)=W_{3,n}^{G}(t)-W_{4,n}^{G}(t)$, with
$W_{3,n}^{G}(t)=n^{-1/2}\sum_{j=1}^n\{\exp(t^\top\varepsilon_j)-m(t)\}$ and
$W_{4,n}^{G}(t)=m(t)a(t)^\top n^{-1/2} \sum_{j=1}^nL_j$. The mean value theorem gives
\[\mathbb E \left[\{\exp(t^\top \varepsilon)-\exp(s^\top \varepsilon)\}^2\right] \leq \kappa\|t-s\|^2, \quad s,t \in B(K),
\]
for some positive $\kappa$. From Theorem 12.3 in Billingsley (1968), the process $\{W_{3,n}^G(t), t \in B(K) \}$  is tight.
By the central limit theorem for martingale differences, $n^{-1/2} \sum_{j=1}^nL_j$ converges in law to a $v$-variate zero mean normal random vector.
 Hence $\{W_{4,n}^G(t), t \in B(K) \}$, being a product of a continuous function and a term which is $O_{\mathbb{P}}(1)$, is tight, and the same property holds for
 $\{W_{2,n}^G(t), t \in B(K) \}$.

\noindent \underline{Proof of (b)}:  In view of
$\mathbb E \left[W_{2,n}^G(t)^2 \right] = \mathbb E \left[U_1(t)^2\right]<\infty$, for each $\varepsilon>0$ there is some positive constant $K$ so that \eqref{resto1} holds.
Likewise,  \eqref{resto2} holds, which completes the proof.  $\bewend$

\medskip

\noindent {\sc Proof} of Theorem \ref{thmconsistG}.  Let  $\varepsilon_j(\theta)=\Sigma_j^{-1/2}(\theta)X_j$. Notice that  $\varepsilon_j(\theta_0)=\varepsilon_j$. Let
$\widetilde{M}_n^G(t)= n^{-1} \!\sum_{j=1}^n \! \exp\{t^\top \widetilde{\varepsilon}_j(\hat{\theta}_n) \}$,
$\widehat{M}_n^G(t)=n^{-1}\! \sum_{j=1}^n \! \exp\{t^\top \varepsilon_j(\hat{\theta}_n) \}$,
$M_n^\circ(t)=n^{-1}\! \sum_{j=1}^n \! \exp\{t^\top \varepsilon_j \}$ and
$B(K):=\{t \in \RR^d:\|t\| \le K\}$. To show the result we will prove
\begin{itemize} \itemsep=0pt
\item[(a)] $\sup_{t \in B(K)} |\widehat{M}_n^G(t)-M_n^\circ(t)|=o_{\mathbb{P}}(1)$,
\item[(b)] $\sup_{t \in B(K)} |\widetilde{M}_n^G(t)-\widehat{M}_n^G(t)|=o_{\mathbb{P}}(1)$,
\end{itemize}
and the result will follow by using the same proof as in the i.i.d. case.

\noindent \underline{Proof of (a)}: Let $\widehat{\theta}_n\! =\! (\widehat{\theta}_{n1}, \ldots, \widehat{\theta}_{nv})^\top$,
$\theta_0\! = \! (\theta_{01}, \ldots, \theta_{0v})^\top$ and $A_{jk}(\theta)\! =\! \Sigma_j^{-1/2}\! (\theta)\frac{\partial}{\partial \theta_k} \Sigma_j^{1/2}\! (\theta)$.
We have  $\varepsilon_j(\widehat{\theta}_n) =\varepsilon_j+\Delta_{n,j}$, with $\Delta_{n,j}=-\sum_{k=1}^vA_{jk}(\widetilde{\theta}_{n,j})\varepsilon_j(\widehat{\theta}_{nk}-\theta_{0k})$,
 for some $\widetilde{\theta}_{n,j}$ between $\widehat{\theta}_{n}$ and $\theta_0$.  Observe that
 $\exp(t^\top \Delta_{n,j})-1=t^\top \Delta_{n,j}\exp(\alpha_{n,j}t^\top \Delta_{n,j})$
 for some $\alpha_{n,j}\in (0,1)$. Now (A.1) and (A.6) yield
 $\| \Delta_{n,j}\| \leq D_j \|\varepsilon_j\|\| \widehat{\theta}_{n}-\theta_0\|$ for large enough $n$, where
 $ \mathbb{E} (D_j^2)<\infty$. The Cauchy--Schwarz inequality gives %We have that
 \[ |\widehat{M}_n^G(t)-M_n^\circ(t)|=\left|\frac{1}{n}\sum_{j=1}^n \exp(t^\top \varepsilon_j)\left\{\exp(t^\top \Delta_{n,}j)-1\right\}\right| \leq
 r_{1,n}(t)^{1/2}r_{2,n}(t)^{1/2},\]
 where $r_{1,n}(t)=M_n(2t)$, and
 \[r_{2,n}(t)=\|t\|^2 \|\widehat{\theta}_n-\theta_0\|^2\exp\left\{2\|t\| \|\widehat{\theta}_n-\theta_0\| \max_{1\leq j \leq n}D_j \|\varepsilon_j\|
 \right\} \frac{1}{n}\sum_{j=1}^n D_j^2 \|\varepsilon_j\|^2.\]
 From the strong law of large numbers in the  Banach space of continuous functions on $B(K)$, we have
 \[
 \sup_{t\in B(K)} r_{1,n}(t) \leq \sup_{t\in B(K)} M_\varepsilon(2t) +\sup_{t\in B(K)}|M_n^G(2t)+M_\varepsilon(2t)| <K_1 \quad \PP\textrm{-a.s.}
 \]
 for some positive constant $K_1$.
 From the ergodic theorem, $n^{-1}\sum_{j=1}^n D_j^2 \|\varepsilon_j\|^2<K_2$ $\PP$-almost surely for some positive constant $K_2$.
 Using stationarity and finite second moments, if follows that
$ \max_{1\leq j \leq n}D_j \|\varepsilon_j\|/\sqrt{n} \to 0$, $\PP$-almost surely. Hence (A.1) yields
$\sup_{t \in B(K)} r_{2,n}(t) \to 0$, in probability. This concludes the proof of (a).

\noindent \underline{Proof of (b)}: \hspace{2pt} The reasoning follows similar steps as the proof of fact (c.1) in the proof of Theorem 7.1 in Henze et al. (2017) and is thus omitted.  $\bewend$

%%%%%%%%%%%%%%%%%%%%%%%%%%%%%%%%%%%%%%%%%%%%%%%%%%%%%%%%%%%%%%%%%%%%%%%%%%%%%%%%%%%%%%%%%%%%%%%%%%%%%%%%%%%%%%%%%%%%%%%%%%%%%%%%%%%%%%%%%%%%%%%%%%%%%%%%%%%%%%%%%%%%%%%%%%%%%%%%%%%%%%%%%%%%%%%%%%%%%%%
\section{Monte Carlo results}\label{sec_monte}
This section describes and summarizes the results of an extensive simulation experiment carried out to study the finite-sample performance of the proposed tests.
Moreover, we consider a real data set of monthly log returns. All computations have
been performed using programs written in the R language.

\subsection{Numerical experiments for i.i.d. data}
Upper quantiles of the null distribution of $T_{n,\beta}$ have been approximated by  gene\-ra\-ting 100,000 samples  from a law N$_d(0,\textrm{I}_d)$. Table  \ref{critical.points}  displays some critical values with the convention that an entry like $^{-4}1.17$ stands for $1.17 \times 10^{-4}$.
The results show that large sample sizes are required to approximate  the critical values by their
corresponding asymptotic values.

\begin{table} %[htb]
\vspace*{-2cm}
\centering
\caption{Critical points for $\pi^{-d/2}T_{n,\beta}$.} \label{critical.points}
{\small
\begin{tabular}{|ccc|rrrrrrr|}
   \hline
    &     &          & \multicolumn{7}{c|}{$\beta$}\\ \cline{4-10}
$d$ & $n$ & $\alpha$ & \multicolumn{1}{c}{2.5} & \multicolumn{1}{c}{3.0} & \multicolumn{1}{c}{3.5} & \multicolumn{1}{c}{4.0} &
\multicolumn{1}{c}{5.0} &  \multicolumn{1}{c}{6.0} &  \multicolumn{1}{c|}{10.0}  \\ \hline
2  & 20 & 0.05 &  0.213 & $^{-1}$0.751 &  $^{-2}$3.269 & $^{-2}$1.639 &  $^{-3}$5.408 & $^{-3}$2.266 & $^{-4}$2.241  \vspace{-3pt}\\
   &    & 0.10 &  0.339 & $^{-1}$1.147 &  $^{-2}$4.857 & $^{-2}$2.380 &  $^{-3}$7.638 & $^{-3}$3.150 & $^{-4}$3.025  \vspace{-3pt}\\
   & 50 & 0.05 &  0.391 & $^{-1}$1.246 &  $^{-2}$5.098 & $^{-2}$2.436 &  $^{-3}$7.594 & $^{-3}$3.078 & $^{-4}$2.875  \vspace{-3pt}\\
   &    & 0.10 &  0.661 & $^{-1}$1.997 &  $^{-2}$7.802 & $^{-2}$3.624 & $^{-3}$10.917 & $^{-3}$4.330 & $^{-4}$3.897  \vspace{-3pt}\\
  & 100 & 0.05 &  0.511 & $^{-1}$1.539 &  $^{-2}$6.073 & $^{-2}$2.838 &  $^{-3}$8.620 & $^{-3}$3.429 & $^{-4}$3.111  \vspace{-3pt}\\
   &    & 0.10 &  0.868 & $^{-1}$2.432 &  $^{-2}$9.168 & $^{-2}$4.153 & $^{-3}$12.094 & $^{-3}$4.724 & $^{-4}$4.143  \vspace{-3pt}\\
  & 200 & 0.05 &  0.612 & $^{-1}$1.757 &  $^{-2}$6.719 & $^{-2}$3.085 &  $^{-3}$9.181 & $^{-3}$3.616 & $^{-4}$3.232  \vspace{-3pt}\\
  &     & 0.10 &  1.028 & $^{-1}$2.726 &  $^{-2}$9.908 & $^{-2}$4.382 & $^{-3}$12.528 & $^{-3}$4.845 & $^{-4}$4.221  \vspace{-3pt}\\
  & 300 & 0.05 &  0.679 & $^{-1}$1.894 &  $^{-2}$7.114 & $^{-2}$3.223 &  $^{-3}$9.466 & $^{-3}$3.719 & $^{-4}$3.296  \vspace{-3pt}\\
  &     & 0.10 &  1.132 & $^{-1}$2.878 & $^{-2}$10.259 & $^{-2}$4.518 & $^{-3}$12.748 & $^{-3}$4.905 & $^{-4}$4.232  \vspace{-3pt}\\
  & 400 & 0.05 &  0.701 & $^{-1}$1.925 &  $^{-2}$7.165 & $^{-2}$3.248 &  $^{-3}$9.502 & $^{-3}$3.721 & $^{-4}$3.283  \vspace{-3pt}\\
  &     & 0.10 &  1.148 & $^{-1}$2.868 & $^{-2}$10.084 & $^{-2}$4.417 & $^{-3}$12.521 & $^{-3}$4.843 & $^{-4}$4.187  \vspace{-3pt}\\
%  & 500 & 0.05 &  0.723 & $^{-1}$1.958 &  $^{-2}$7.276 & $^{-2}$3.288 &  $^{-3}$9.568 & $^{-3}$3.737 & $^{-4}$3.310  \vspace{-3pt}\\
%  &     & 0.10 &  1.171 & $^{-1}$2.903 & $^{-2}$10.217 & $^{-2}$4.451 & $^{-3}$12.577 & $^{-3}$4.846 & $^{-4}$4.200  \vspace{-3pt}\\
\hline
3  & 20 & 0.05 &  0.356 & $^{-1}$1.066 &  $^{-2}$4.095 & $^{-2}$1.851 &  $^{-3}$5.218 & $^{-3}$1.942 & $^{-4}$1.413   \vspace{-3pt}\\
   &    & 0.10 &  0.520 & $^{-1}$1.504 &  $^{-2}$5.629 & $^{-2}$2.503 &  $^{-3}$6.886 & $^{-3}$2.518 & $^{-4}$1.773   \vspace{-3pt}\\
   & 50 & 0.05 &  0.719 & $^{-1}$1.906 &  $^{-2}$6.760 & $^{-2}$2.894 &  $^{-3}$7.598 & $^{-3}$2.709 & $^{-4}$1.828   \vspace{-3pt}\\
   &    & 0.10 &  1.153 & $^{-1}$2.879 &  $^{-2}$9.789 & $^{-2}$4.073 & $^{-3}$10.317 & $^{-3}$3.593 & $^{-4}$2.334    \vspace{-3pt}\\
  & 100 & 0.05 &  0.988 & $^{-1}$2.433 &  $^{-2}$8.258 & $^{-2}$3.426 &  $^{-3}$8.696 & $^{-3}$3.043 & $^{-4}$1.992   \vspace{-3pt}\\
   &    & 0.10 &  1.646 & $^{-1}$3.732 & $^{-2}$11.943 & $^{-2}$4.788 & $^{-3}$11.572 & $^{-3}$3.945 & $^{-4}$2.489    \vspace{-3pt}\\
  & 200 & 0.05 &  1.232 & $^{-1}$2.851 &  $^{-2}$9.322 & $^{-2}$3.781 &  $^{-3}$9.365 & $^{-3}$3.231 & $^{-4}$2.078    \vspace{-3pt}\\
  &     & 0.10 &  2.046 & $^{-1}$4.319 & $^{-2}$13.243 & $^{-2}$5.167 & $^{-3}$12.210 & $^{-3}$4.123 & $^{-4}$2.567     \vspace{-3pt}\\
  & 300 & 0.05 &  1.332 & $^{-1}$2.979 &  $^{-2}$9.555 & $^{-2}$3.849 &  $^{-3}$9.431 & $^{-3}$3.242 & $^{-4}$2.072    \vspace{-3pt}\\
  &     & 0.10 &  2.187 & $^{-1}$4.441 & $^{-2}$13.364 & $^{-2}$5.156 & $^{-3}$12.105 & $^{-3}$4.073 & $^{-4}$2.527    \vspace{-3pt}\\
  & 400 & 0.05 &  1.397 & $^{-1}$3.061 &  $^{-2}$9.725 & $^{-2}$3.893 &  $^{-3}$9.509 & $^{-3}$3.260 & $^{-4}$2.084     \vspace{-3pt}\\
  &     & 0.10 &  2.245 & $^{-1}$4.481 & $^{-2}$13.341 & $^{-2}$5.122 & $^{-3}$12.010 & $^{-3}$4.046 & $^{-4}$2.519    \vspace{-3pt}\\
%  & 500 & 0.05 &  1.441 & $^{-1}$3.106 &  $^{-2}$9.799 & $^{-2}$3.908 &  $^{-3}$9.497 & $^{-3}$3.254 & $^{-4}$2.073    \vspace{-3pt}\\
%  &     & 0.10 &  2.296 & $^{-1}$4.480 & $^{-2}$13.253 & $^{-2}$5.089 & $^{-3}$11.932 & $^{-3}$4.003 & $^{-4}$2.496   \vspace{-3pt}\\
\hline
5  & 20 & 0.05 &  0.597 & $^{-1}$1.347 &  $^{-2}$4.130 & $^{-2}$1.554 &  $^{-3}$3.275 & $^{-3}$0.971 & $^{-5}$3.884 \vspace{-3pt}\\
   &    & 0.10 &  0.774 & $^{-1}$1.691 &  $^{-2}$5.089 & $^{-2}$1.886 &  $^{-3}$3.900 & $^{-3}$1.142 & $^{-5}$4.474 \vspace{-3pt}\\
   & 50 & 0.05 &  1.519 & $^{-1}$2.868 &  $^{-2}$7.862 & $^{-2}$2.731 &  $^{-3}$5.215 & $^{-3}$1.460 & $^{-5}$5.283 \vspace{-3pt}\\
   &    & 0.10 &  2.332 & $^{-1}$4.130 & $^{-2}$10.801 & $^{-2}$3.633 &  $^{-3}$6.633 & $^{-3}$1.809 & $^{-5}$6.260 \vspace{-3pt}\\
  & 100 & 0.05 &  2.315 & $^{-1}$3.947 & $^{-2}$10.132 & $^{-2}$3.381 &  $^{-3}$6.134 & $^{-3}$1.667 & $^{-5}$5.768 \vspace{-3pt}\\
   &    & 0.10 &  3.782 & $^{-1}$5.884 & $^{-2}$14.199 & $^{-2}$4.530 &  $^{-3}$7.779 & $^{-3}$2.051 & $^{-5}$6.782 \vspace{-3pt}\\
  & 200 & 0.05 &  3.047 & $^{-1}$4.744 & $^{-2}$11.541 & $^{-2}$3.736 &  $^{-3}$6.565 & $^{-3}$1.755 & $^{-5}$5.962 \vspace{-3pt}\\
  &     & 0.10 &  4.969 & $^{-1}$6.964 & $^{-2}$15.880 & $^{-2}$4.896 &  $^{-3}$8.131 & $^{-3}$2.112 & $^{-5}$6.875 \vspace{-3pt}\\
  & 300 & 0.05 &  3.346 & $^{-1}$5.016 & $^{-2}$11.974 & $^{-2}$3.829 &  $^{-3}$6.636 & $^{-3}$1.769 & $^{-5}$5.985 \vspace{-3pt}\\
  &     & 0.10 &  5.445 & $^{-1}$7.343 & $^{-2}$16.307 & $^{-2}$4.960 &  $^{-3}$8.119 & $^{-3}$2.100 & $^{-5}$6.832 \vspace{-3pt}\\
  & 400 & 0.05 &  3.608 & $^{-1}$5.234 & $^{-2}$12.292 & $^{-2}$3.889 &  $^{-3}$6.679 & $^{-3}$1.776 & $^{-5}$5.997 \vspace{-3pt}\\
  &     & 0.10 &  5.838 & $^{-1}$7.586 & $^{-2}$16.477 & $^{-2}$4.958 &  $^{-3}$8.085 & $^{-3}$2.085 & $^{-5}$6.821 \vspace{-3pt}\\
%  & 500 & 0.05 &  3.746 & $^{-1}$5.330 & $^{-2}$12.398 & $^{-2}$3.906 &  $^{-3}$6.690 & $^{-3}$1.776 & $^{-5}$5.999 \vspace{-3pt}\\
%  &     & 0.10 &  6.048 & $^{-1}$7.649 & $^{-2}$16.484 & $^{-2}$4.934 &  $^{-3}$7.999 & $^{-3}$2.066 & $^{-5}$6.781  \vspace{-3pt}\\
\hline
\end{tabular}} \end{table}

A natural competitor of the test based on $T_{n,\beta}$ is the CF-based test studied in  Henze and Wagner (1997) (HW-test).
The latter procedure is simple to compute as well as  affine invariant, and it has revealed good power performance with regard to competitors.
The behaviour of the test based on $T_{n,\beta}$ in relation to the HW-test depends on whether the distribution is heavy-tailed or not. We tried  a number of non-heavy-tailed distributions (specifically, the multivariate Laplace distribution,  finite mixtures of normal distributions, the skew-normal distribution, the multivariate $\chi^2$-distribution,
%the generalized exponential distribution,
the Khintchine distribution, the uniform distribution on $[0,1]^d$ and the Pearson type II family). For these distributions we observed that the power of the proposed test is either similar
or smaller than that of the HW-test; for  very heavy-tailed distributions,  the new test outperforms the HW-test. This observation can be appreciated by looking at
Table \ref{tabla3}, which displays the empirical power calculated by generating 10,000 samples (in each case), for the significance level $\alpha=0.05$,
from the following heavy-tailed alternatives:
($ASE_{\theta}$) the $\theta$-stable  and  elliptically-contoured distribution and the
($T_{\theta}$) multivariate Student's $t$ with $\theta$ degrees of freedom.
The same fact was also observed in Zghoul (2010), who numerically studied  the test based on $T_{n,\beta}$ for univariate data.

In  our simulations we tried a large number of values for $\beta$ %ranging from 2.1 to 20
for the proposed test as well as for the HW-test. % (in the last case $\beta$ was taken between 0.1 to 10).
 The tables display the results for those values of $\beta$ giving the highest power in most of the cases considered. The same can be said for the simulations in the next subsection.

\begin{table} %[htb]
\vspace*{-2cm}
\centering
\caption{Percentage of rejection for nominal level $\alpha=0.05$ and $n=50$.} \label{tabla3}
\footnotesize
\begin{tabular}{|cc|rrrrrr|rrr|}
\hline
     & &   \multicolumn{6}{c|}{Test based on $T_{n,\beta}$} & \multicolumn{3}{c|}{HW-test }\\ \cline{3-11}
     & &   \multicolumn{6}{c|}{$\beta$} & \multicolumn{3}{c|}{$\beta$}\\  \hline

   & $d$ &   3&  3.5 &  4.0 & 5.0 & 6.0 & 10.0 & 0.1 & 0.5 & 1.0\\  \hline
$ASE_{1.75}$
 &  2 &  72.47 & 72.62 & 72.43 & 72.08 & 71.59 & 70.34 & 67.29 & 67.75 & 59.91  \vspace{-3pt}\\
 &  3 &  82.70 & 82.78 & 82.76 & 82.69 & 82.52 & 81.92 & 79.07 & 78.16 & 68.60  \vspace{-3pt}\\
 &  5 &  90.51 & 90.86 & 90.95 & 91.25 & 91.32 & 91.05 & 88.89 & 87.46 & 75.71  \vspace{-3pt}\\
\hline
$ASE_{1.85}$
 & 2 & 54.43 & 54.39 & 54.35 & 53.91 & 53.55 & 52.59 & 50.00 & 48.17 & 39.35  \vspace{-3pt}\\
 & 3 & 62.72 & 62.67 & 62.61 & 62.46 & 62.39 & 61.58 & 57.95 & 54.67 & 42.44  \vspace{-3pt}\\
 & 5 & 75.31 & 75.52 & 75.65 & 75.96 & 76.03 & 75.63 & 71.81 & 66.66 & 47.82  \vspace{-3pt}\\
\hline
$ASE_{1.95}$
 & 2 & 24.67 & 24.62 & 24.52 & 24.22 & 24.11 & 23.56 & 22.44 & 20.78 & 15.38  \vspace{-3pt}\\
 & 3 & 29.31 & 29.37 & 29.47 & 29.12 & 28.91 & 28.45 & 26.37 & 24.04 & 16.79  \vspace{-3pt}\\
 & 5 & 38.28 & 38.39 & 38.37 & 38.27 & 38.00 & 37.55 & 33.99 & 29.35 & 17.39  \vspace{-3pt}\\
\hline
$T_{5}$
 & 2 & 58.77 & 58.82 & 58.74 & 58.26 & 57.82 & 56.21 & 51.44 & 54.20 & 47.58  \vspace{-3pt}\\
 & 3 & 40.59 & 40.76 & 40.98 & 41.11 & 41.34 & 41.31 & 39.69 & 37.79 & 28.99  \vspace{-3pt}\\
 & 5 & 87.14 & 87.80 & 88.43 & 89.17 & 89.59 & 89.76 & 86.36 & 87.21 & 77.92  \vspace{-3pt}\\
\hline
$T_{7}$
 & 2 & 42.33 & 42.21 & 42.16 & 41.86 & 41.43 & 39.68 & 36.19 & 36.30 & 28.97  \vspace{-3pt}\\
 & 3 & 55.23 & 55.49 & 55.51 & 55.62 & 55.30 & 54.38 & 49.20 & 48.82 & 37.89  \vspace{-3pt}\\
 & 5 & 71.51 & 72.50 & 73.17 & 74.07 & 74.47 & 74.30 & 69.35 & 67.75 & 51.33  \vspace{-3pt}\\
\hline
$T_{10}$
 & 2 & 28.80 & 28.82 & 28.73 & 28.29 & 27.99 & 27.12 & 24.48 & 22.94 & 16.23  \vspace{-3pt}\\
 & 3 & 38.34 & 38.56 & 38.58 & 38.51 & 38.33 & 37.00 & 32.97 & 30.55 & 20.84  \vspace{-3pt}\\
 & 5 & 51.64 & 52.38 & 53.01 & 53.91 & 54.25 & 54.33 & 48.85 & 45.36 & 28.41  \\
\hline
\end{tabular}
\end{table}

\subsection{Numerical experiments for GARCH data}
In our simulations we considered a bivariate CCC--GARCH(1,1) model with
\[ {b}=\left(\begin{array}{c}0.1\\ 0.1\end{array}\right), \
 {{B}}_1= \left(\begin{array}{cc}0.1 & 0.1\\ 0.1 & 0.1\end{array}\right), \
 {{\Gamma}}_1= \left(\begin{array}{cc} \gamma & 0.01\\ 0.01 & \gamma \end{array}\right), \
 {{R}}=\left(\begin{array}{cc}1 & r\\ r & 1 \end{array}\right),
\]
for $\gamma=0.3, 0.4, 0.5$ and $r=0, 0.3$, and a trivariate CCC--GARCH(1,1) model with  ${b}=(0.1, 0.1, 0.1)'$,
\[
 {B}_1= \left(\begin{array}{ccc}0.1 & 0.1 & 0.1\\ 0.1 & 0.1 & 0.1 \\ 0.1 & 0.1 & 0.1\end{array}\right), \
 {\Gamma}_1= \left(\begin{array}{ccc} \gamma & 0.01 & 0.01 \\ 0.01 & \gamma & 0.01\\  0.01 & 0.01 & \gamma \end{array}\right), \
 {R}=\left(\begin{array}{ccc}1 & r & r\\ r & 1 & r\\ r & r & 1 \end{array}\right)
\]
and $\gamma$ and $r$ as before. The parameters in the CCC-GARCH models were estimated by their QMLE using the package {\tt ccgarch} of the language R. For the distribution of the innovations, we took  ${\varepsilon}_1, \ldots, {\varepsilon}_n$ i.i.d. from the distribution of ${\varepsilon}$ with
${\varepsilon}$ having a ($N$)  multivariate normal distribution, in order to study the level of the resulting bootstrap test. To assess
 the power we considered the following heavy-tailed distributions:  $T_\theta$, the multivariate $\beta$-generalized distribution ($GN_{\theta}$), that coincides with the normal distribution for $\theta=2$ and has heavy tails for $0<\theta<2$
(Goodman and Kotz, 1973),  and
the asymmetric exponential power distribution ($AEP$), whereby $(X_1, \ldots, X_d)^\top$, with $X_1,\ldots,  X_d$ i.i.d. from a univariate $AEP$ distribution  (Zhu and Zinde-Walsh, 2009) with parameters $\alpha=0.4$,  $p_1=1.182$ and $p_2=1.820$ (these settings gave useful results in practical applications for the  errors in GARCH type models). As in the previous subsection, we also calculated the HW-test.

Table \ref{tablagarch}  reports the percentages of rejections for nominal significance level $\alpha=0.05$ and sample size $n=300$, for $r=0,\, 0.3$ and $\gamma=0.4$. The resulting pictures for $\gamma=0.3, \, 0.5$ are quite similar so, to save space, we omit the results for these values of $\gamma$.
In order to reduce the computational burden we adopted the warp-speed method of  Giacomini et al. (2013), which works as follows:
rather than computing critical points for each Monte Carlo sample, one resample is generated for each Monte Carlo sample,
and the resampling test statistic is computed for that sample; then the resampling critical values for $T^G_{n,\beta}$ are computed from the empirical distribution
determined by the resampling replications of $T_{n,\beta}^{G*}$. In our simulations we generated $10,000$ Monte Carlo samples for the level and $2,000$ for the power.
Looking at Table  \ref{tablagarch}, we conclude that: the actual level of the proposed bootstrap test is very close to the nominal level, and this is also true for the HW-test (although to the best of our knowledge,  the consistency of the bootstrap null distribution estimator of the HW-test statistic has been proved only for the univariate case in Jim\'enez-Gamero, 2014); and with respect to the power, the proposed test in most cases outperforms the HW-test.

\begin{table} %[htb]
%\vspace*{-2cm}
\centering
\caption{Percentage of rejections for nominal level $\alpha=0.05$, $\gamma=0.4$ and $n=300$.} \label{tablagarch}
\footnotesize
\begin{tabular}{|ccc|rrrrr|rrrr|}
\hline
 &   &  &  \multicolumn{5}{c|}{Test based on $T^G_{n,\beta}$} & \multicolumn{4}{c|}{HW-test }\\ \cline{4-12}
 &   &  &   \multicolumn{5}{c|}{$\beta$} & \multicolumn{4}{c|}{$\beta$}\\  \hline
    & $d$      & $r$   &  2.1 &  2.2 &  2.3 &  2.4 & 2.5  &  1.0 & 1.5 & 2.0 & 2.5 \\ \hline
$N$ &  2       & 0.0   & 4.96 & 4.85 & 4.81 & 4.79 & 4.73 & 5.06 & 4.80 & 4.97 & 4.82   \vspace{-3pt}\\
    &          & 0.3   & 4.14 & 4.33 & 4.38 & 4.40 & 4.27 &  4.95 & 5.45 & 5.36 & 5.29\vspace{-3pt} \\
    &  3       & 0.0   & 4.54 & 4.71 & 4.73 & 4.74 & 4.73 &  4.64 & 4.64 & 4.88 & 4.51 \vspace{-3pt}\\
    &          & 0.3   & 4.96 & 4.85 & 4.81 & 4.79 & 4.73 &  5.06 & 4.80 & 4.97 & 4.82  \vspace{-3pt}\\ \hline
$T_{10}$ & 2   & 0.0   & 61.85 & 61.20 & 59.25 & 57.55 & 55.50  & 26.70 & 36.70 & 37.20 & 34.85    \vspace{-3pt} \\
    &          & 0.3   & 66.95 & 66.80 & 65.85 & 64.15 & 61.35  & 20.50 & 31.70 & 32.10 & 30.60   \vspace{-3pt} \\
    &    3     & 0.0   & 81.45 & 80.95 & 80.15 & 79.65 & 78.20  & 45.75 & 55.40 & 50.95 & 43.80        \vspace{-3pt} \\
    &          & 0.3   & 78.30 & 78.05 & 78.20 & 77.20 & 77.15  & 42.40 & 55.70 & 52.85 & 44.00  \vspace{-3pt}\\ \hline
$GN_{1.65}$ & 2 & 0.0  & 22.40 & 21.05 & 20.10 & 18.95 & 17.85  &  8.65 & 15.20 & 16.45 & 16.75  \vspace{-3pt} \\
    &          & 0.3   & 18.30 & 17.80 & 17.70 & 16.80 & 16.10  &  8.00 & 14.00 & 16.00 & 14.30  \vspace{-3pt} \\
    &  3       & 0.0   & 17.55 & 18.40 & 18.10 & 17.80 & 16.90  &  9.10 & 14.85 & 15.35 & 15.60  \vspace{-3pt} \\
    &          & 0.3   & 20.00 & 19.65 & 19.85 & 19.80 & 18.90  &  9.70 & 13.95 & 15.55 & 15.15  \vspace{-3pt} \\  \hline
$AEP$      & 2 & 0.0   & 56.75 & 55.50 & 53.35 & 51.10 & 49.00  & 29.55 & 49.85 & 52.85 & 51.45  \vspace{-3pt} \\
    &          & 0.3   & 52.70 & 51.20 & 49.65 & 47.90 & 45.75  & 26.35 & 45.20 & 50.00 & 49.20  \vspace{-3pt} \\
    &     3    & 0.0   & 55.40 & 55.85 & 54.85 & 53.75 & 51.65  & 38.25 & 54.25 & 55.45 & 49.25  \vspace{-3pt} \\
    &          & 0.3   & 59.55 & 59.30 & 58.75 & 57.15 & 57.00  & 33.15 & 53.65 & 53.90 & 49.70  \vspace{-3pt}\\   \hline
\end{tabular}
\end{table}

\subsection{A real data set application}
\begin{figure}
\begin{center}
\includegraphics[width=7.9cm, height=8.8cm]{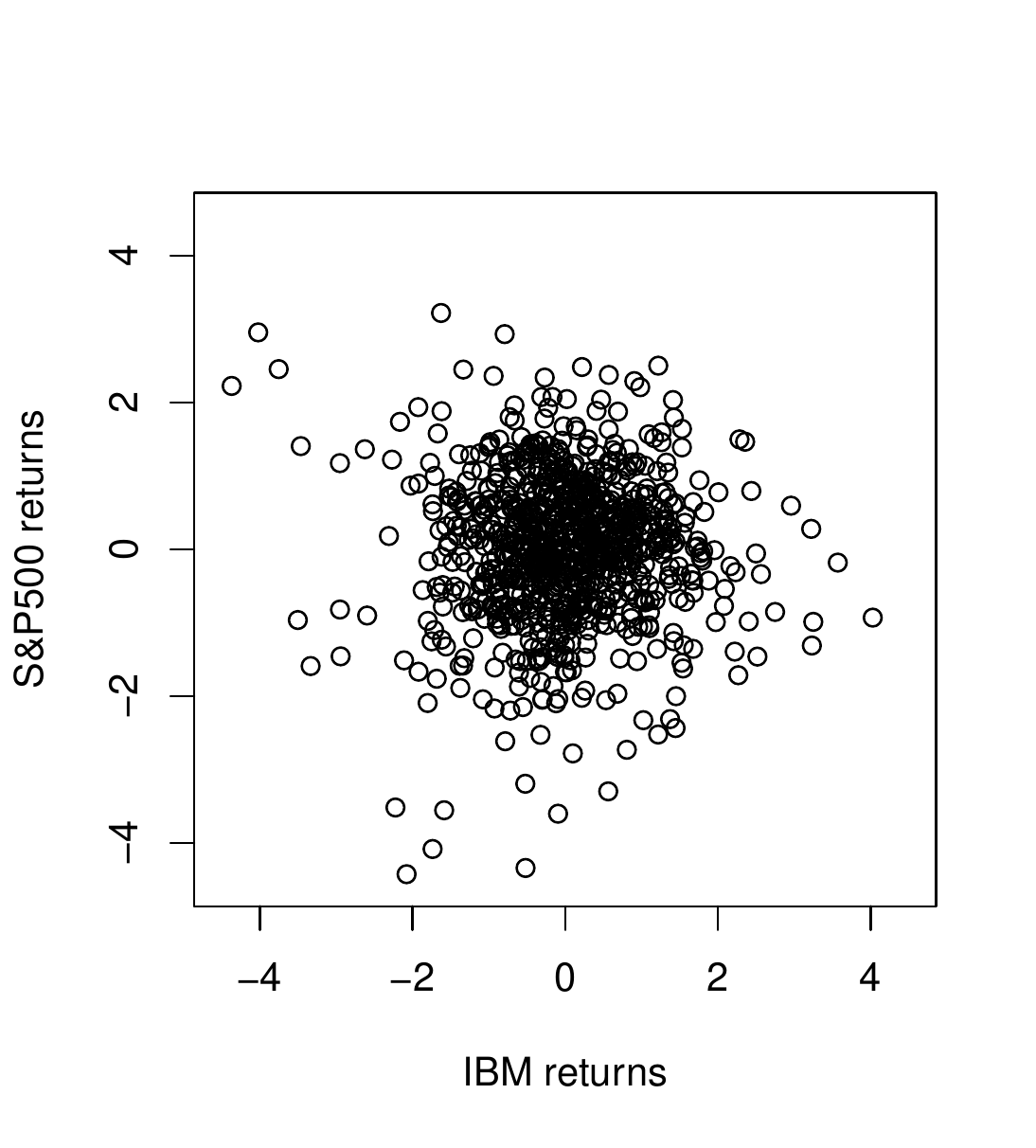}
\end{center}
\caption{Scatter plot of the residuals.} \label{fig1}
\end{figure}

\begin{figure}
\begin{center}
\includegraphics[width=15.2cm, height=7.1cm]{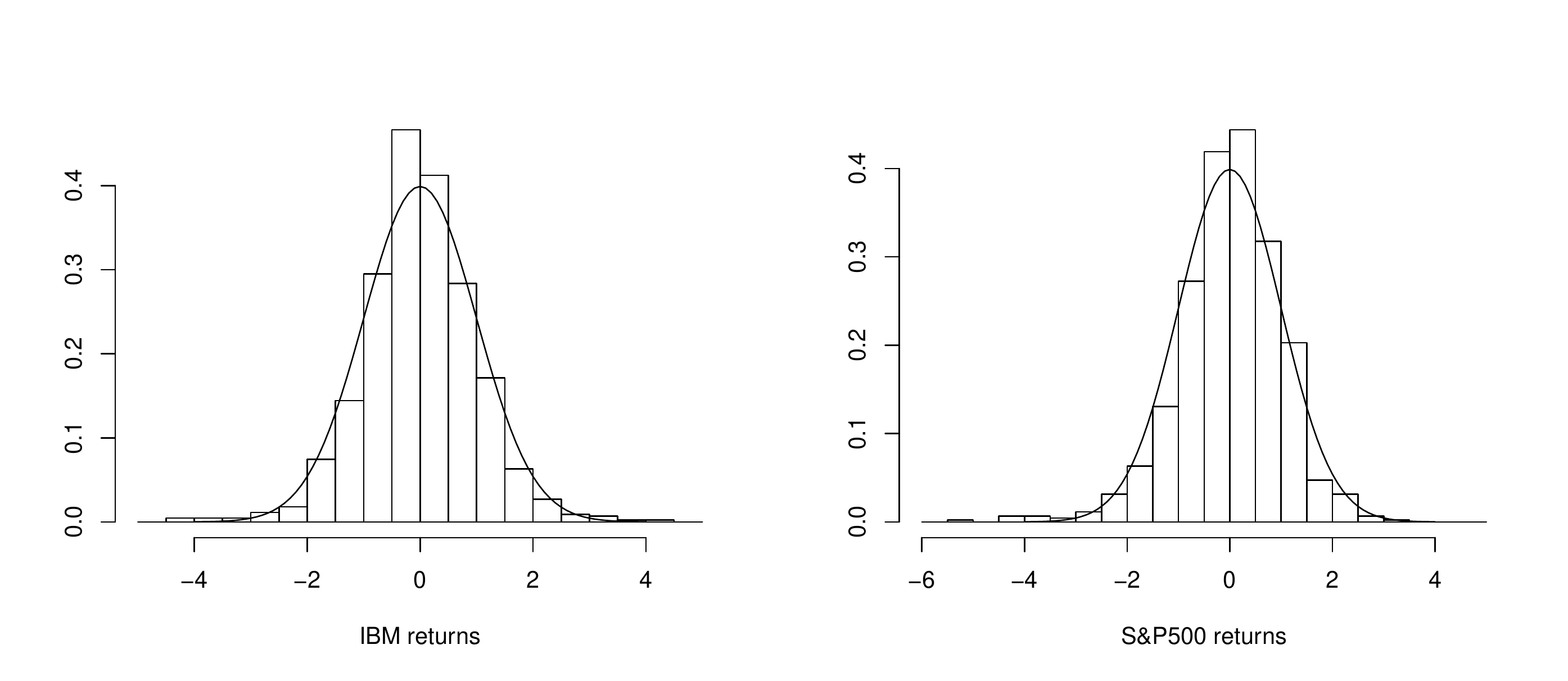}
\caption{Histograms  of the residuals.} \label{fig2}
\end{center}
\end{figure}

As an illustration, we consider
the monthly log returns of IBM stock and the S\&P 500
index from January 1926 to December 2008 with 888 observations. This data set was analyzed in Example 10.5 of Tsay(2010), where it is showed that a CCC-GARCH(1,1) model provides a adequate description of the data,  which is available from the  website  {\tt http://faculty.chicagobooth.edu/ruey.tsay/teaching/fts/}) of the author. We applied the proposed test and the HW test for testing $H_{0,G}$. The $p$-values were obtained by generating 1000 bootstrap samples. For all values of $\beta$ in Table  \ref{tablagarch} we get the same $p$-value, 0.000, which leads us to reject $H_{0,G}$, as expected by looking at Figure \ref{fig1}, which displays the scatter plot of the residuals after fitting a  CCC-GARCH(1,1) model to the log returns, and Figure \ref{fig2}, that  represents  the histograms of the marginal residuals with the probability density function of a standard normal law superimposed.

%%%%%%%%%%%%%%%%%%%%%%%%%%%%%%%%%%%%%%%%%%%%%%%%%%%%%%%%%%%%%%%%%%%%%%%%%%%%%%%%%%%%%%%%%%%%%%%%%%%%%%%%%%%%%%%%%%%%%%%%%%%%%%%%%%%%%%%%%%%%%%%%%%%%%%%%%%%%%%%%%%%%%%%%%%%%%%%%%%
%
%
%
%
%
%%%%%%%%%%%%%%%%%%%%%%%%%%%%%%%%%%%%%%%%%%%%%%%%%%%%%%%%%%%%%%%%%%%%%%%%%%%%%%%%%%%%%%%%%%%%%%%%%%%%%%%%%%%%%%%%%%%%%%%%%%%%%%%%%%%%%%%%%%%%%%%%%%%%%%%%%%%%%%%%%%%%%%%%%%%%%%%%%%%%%%%%%%%

\section{Conclusions}\label{secconclusion}
We have studied a class of affine invariant  tests for multivariate normality both in an i.i.d. setting and in the context of
testing that the innovation distribution of a multivariate GARCH model is Gaussian, thus generalizing results of Henze and Koch (2017) in
two ways. The test statistics are suitably weighted $L^2$-statistics
based on the difference between the empirical moment generating function of scaled residuals of the data and the moment generating function of the standard
normal distribution in $\RR^d$. As such, they can be considered as 'moment generating function analogues'  to the time-honored class of BHEP tests
that use the empirical {\em characteristic} function. As the decay of a weight function figuring in the  test statistic tends to infinity, the test
statistic approaches a certain linear combination of two well-known measures of multivariate skewness.
The tests are easy to implement,  and they turn out to be consistent against a wide range of alternatives. In contrast to a recently studied $L^2$-statistic
of Henze et al. (2017) that uses both the empirical moment generating and the empirical characteristic function, our test is also feasible for larger sample sizes since the
computational complexity is of order $O(n^2)$. Regarding power, the new tests outperform the BHEP-tests against heavy-tailed distributions.

\section*{Acknowledgements}
M.D. Jim\'enez-Gamero was partially supported by
grant  MTM2014-55966-P of the Spanish Ministry of Economy and Competitiveness.

\end{document}